\documentclass[11pt]{amsart}
\usepackage{amssymb}

\usepackage{graphicx}
\usepackage{xcolor}

\newcommand{\binomsq}[2]{\genfrac{[}{]}{0pt}{}{#1}{#2}}

\newcommand{\res}{\upharpoonright}

\theoremstyle{plain}
\newtheorem{theorem}{Theorem}[section]

\newtheorem{lemma}[theorem]{Lemma}
\newtheorem{proposition}[theorem]{Proposition}
\newtheorem{definition}[theorem]{Definition}
\newtheorem{illustration}[theorem]{Illustration}

\theoremstyle{remark}

\numberwithin{equation}{section}

\newcommand{\dbullet}{\,\raisebox{-1.4pt}{\( \bullet \)}\,}
\newcommand{\dowd}{\, .\,}

\begin{document}
\title[Abstract Ramsey theory and Ramsey theorems for trees]{Abstract approach to Ramsey theory\\ and Ramsey theorems for finite trees}

\author{S{\l}awomir Solecki}

\thanks{Research supported by NSF grant DMS-1001623.}

\address{Department of Mathematics\\
University of Illinois\\
1409 W. Green St.\\
Urbana, IL 61801, USA\\
{\rm and }Institute of Mathematics\\
Polish Academy of Sciences\\
ul. {\'S}niadeckich 8\\
00-950 Warsaw, Poland}

\email{ssolecki@math.uiuc.edu}



\begin{abstract} I will give a presentation of an abstract
approach to finite Ramsey theory found in an earlier paper of mine.
I will prove from it a common generalization of Deuber's Ramsey theorem
for regular trees and a recent Ramsey theorem of Jasi{\'n}ski for
boron tree structures. This generalization appears to be new. I will also
show, in exercises, how to deduce from it the Milliken Ramsey
theorem for strong subtrees.
\end{abstract}

\maketitle

\section{Introduction}

The first result of 
pure finite Ramsey theory and a prototype of the many later results of this area (see \cite{Nes}) is 
the theorem proved by Ramsey in 1930. We recall it now to remind the reader of the flavor of pure finite 
Ramsey theory.  
We will also refer to this statement later on. 
For a natural number $n$, let $[n] = \{ 1,\dots, n\}$; in particular, $[0]=\emptyset$. The classical Ramsey 
theorem says that given natural numbers 
$d$, $l$, and $m$, there exists a natural number $n$ such that for each $d$-coloring, that is, a coloring 
with $d$ colors, of all $l$ element subsets of $[n]$, there exists an $m$ element subset $z$ of $[n]$ such 
that all $l$ element subsets of $z$ have the same color.

In Section~\ref{S:aR}, we give an exposition of the abstract approach to pure finite Ramsey
theory developed in \cite{Sol}; the main theorem, saying that a
general pigeonhole principle implies a general Ramsey property, is
stated as Theorem~\ref{T:ram} (see also Appendix~1). 
Most pure finite Ramsey theoretic results can be viewed as instances of the 
machinery presented here.
In the exposition, we make an effort to motivate
the main abstract notions and we also illustrate them with examples. 

In Sections~\ref{S:tree}, \ref{S:JD}, \ref{S:M}, using
arguments consisting mostly of formulating appropriate
definitions, we show that certain Ramsey-theoretic results 
for finite trees, one of which is new,  
are particular instances of the general Theorem~\ref{T:ram}. These
applications of Theorem~\ref{T:ram} to concrete situations are similar
to each other, with the main differences lying in the derivations used (more
on it in the next paragraph). Therefore, in the first two applications
(the illustrations in Section~\ref{S:aR} and
Illustration~\ref{I:tree}), we explicitly check all the details
and provide pictures; in the third application (a generalization of
Deuber's and Jasi{\'n}ski's theorems, Section~\ref{S:JD}), we give
all the definitions, but carefully check the pigeonhole principle
only; in the last application (Milliken's theorem,
Section~\ref{S:M}), we state all the definitions, but leave checking
the details to the reader in exercises. Recall that
in \cite{Sol}, it is shown how, for example, the classical Ramsey theorem,
the Graham--Rothschild theorem,
and a new self-dual Ramsey theorem can be obtained as instances of Theorem~\ref{T:ram}.

In each of the many concrete Ramsey theorems (considered
here and in \cite{Sol}), the same underlying algebraic structure
turns out to be present, the structure of a {\em normed background} given by
Definition~\ref{D:norback} (see also Appendix 1). A crucial element
of such structures is a {\em truncation operator}, which forms a
basis for inductive arguments. In the concrete situations involving
trees and considered in the present paper, there is a close
connection between truncation operators and {\em derivations on
trees}. Roughly speaking there are two natural derivations on trees: cutting
off the rightmost branch and cutting off the highest leaves. These
two derivations give rise to two types of truncation operators,
which lead to two types of normed backgrounds, which in turn lead to
two Ramsey theorems. Namely, the branch cutting derivation gives a
generalization of Deuber's and Jasi{\'n}ski's theorems, while the leaf cutting
derivation gives Milliken's theorem.

{\em For convenience, we adopt the following modification to the
notation for the operation of subtracting $1$ among natural numbers:
we set $0-1$ to be equal to $0$; for $k>0$, $k- 1$ retains its usual
meaning.}

\section{Abstract approach with illustrations}\label{S:aR}

\subsection{Abstract Ramsey theory}

A typical Ramsey-type theorem has the following form. We start with two families $\mathcal F$ and 
$\mathcal P$. (Elements of $\mathcal F$
and $\mathcal P$ are usually finite sets of functions, most frequently some type of morphisms.)
A set $P$ from $\mathcal P$ and a number of colors $d$ are given. The conclusion of the theorem then 
asserts that there is a set $F$ from
$\mathcal F$ with a given mapping (usually a type of composition) defined on $F\times P$, 
\[
F\times P \ni (f,x)\to f\dowd x,
\]
such that for each $d$-coloring of the image $\{ f\dowd
x\colon f\in F,\, x\in P\}$ of $F\times P$ under the mapping there exists $f_0\in F$ with $\{ f_0\dowd
x\colon x\in P\}$ monochromatic. 

Below in the paper, we formalize this vague idea and we also give several concrete examples that should 
convince the 
reader that Ramsey-type theorems do indeed have this form. Here, as an illustration, we only 
phrase the classical Ramsey theorem in a way that is compatible with the general framework above. It 
may be useful for the reader to recall here the Ramsey theorem from the first paragraph of the introduction.
In the restatement of the Ramsey theorem to which we now proceed, for natural numbers $p$ and $q$, 
we identify $p$ element 
subsets of $[q]$ with increasing injections from $[p]$ to $[q]$ so 
that a subset $z$ is identified with the unique increasing injection 
whose range is equal to $z$. One can take ${\mathcal F}={\mathcal P}$ to be the family 
of all sets produced as follows: fix natural numbers $p$ and $q$ and form the set 
of all increasing injections from $[p]$ to $[q]$. 
Fix natural numbers $d$, $l$, and $m$, and let $P\in {\mathcal P}$ be the set of all increasing injections from 
$[l]$ to $[m]$. 
Then the classical Ramsey theorem says that there is an $n$ with the following property. For the set 
$F\in {\mathcal F}$ of all increasing injections 
from $[m]$ to $[n]$, if we $d$-color the set  
\[
\{ f\circ x\colon f\in F,\, x\in P\} = \hbox{ all increasing injections from $[l]$ to $[n]$}, 
\]
then there exists $f_0\in F$ such that $\{ f_0\circ x\colon x\in P\}$ is 
monochromatic. 

Now we start the description of the abstract approach. Let $A$ and $X$ be sets. Assume we are given a {\em partial} function
from $A\times X$ to $X$:
\[
(a,x)\to a\dowd x.
\]
Such a function $\dowd$ will be called an {\bf action} ({\bf of $A$ on
$X$}). No properties of the function $\dowd$ are assumed to hold
at this point. For $F\subseteq A$ and $P\subseteq X$, we say that
$F\dowd P$ {\em is defined} if $a\dowd x$ is defined for all $a\in F$
and $x\in P$, and we let
\[
F\dowd P = \{ a\dowd x\colon a\in F,\, x\in P\}.
\]
We also write $a\dowd P$ for $\{ a\}\dowd P$.

We will give a sequence of illustrations that contain the most
rudimentary examples of the general notions being introduced. The
illustrations depend on each other and lead to the classical Ramsey
theorem.

\begin{illustration}\label{I:ramfir} {\rm Let $A = X$ be the set of
all {\em (strictly) increasing} functions from $[k] = \{ 1, \dots , k\}$ to ${\mathbb N}\setminus \{ 0\}$, where $k$ ranges
over $\mathbb N$. Given $a, x\in A=X$ with $a\colon [l]\to {\mathbb
N}\setminus \{ 0\}$ and $x\colon [k]\to {\mathbb N}\setminus \{ 0\}$, let $a\dowd x$ be defined
precisely when $[l]$ contains the image of $x$ and put
\[
a\dowd x = a\circ x.
\]}
\end{illustration}

Going back to the general situation, let ${\mathcal F}$ and ${\mathcal P}$ be families of 
non-empty subsets of $A$ and $X$, respectively. Assume we have a
partial function from ${\mathcal F}\times {\mathcal P}$ to
${\mathcal P}$:
\[
(F, P)\to F\dbullet P
\]
such that if $F\dbullet P$ is defined, then it is given by the point-wise action
of $F$ on $P$, that is,  $F\dowd P$ is defined
and
\[
F\dbullet P = F\dowd P.
\]
In such a situation, we say that $({\mathcal F}, {\mathcal P},
\dbullet)$ is a {\bf pair of families over $(A,X, \dowd)$}. Introducing
a restriction $\dbullet$ of the point operation of sets in $\mathcal F$ on sets in $\mathcal P$
makes the Ramsey condition (R) below more flexible, while a careful calibration of the resteriction makes it
possible to satisfy condition ($*$) of the next subsection. In concrete situations, definitions of restrictions 
$\dbullet$ are very natural.

\begin{illustration}\label{I:dbul}{\rm For $k,l\in {\mathbb N}$ with $0<k\leq l$, let
$\binom{l}{k}$ stand for the set of all (strictly) increasing functions from
$[k]$ to $[l]$. Let also $\binom{0}{0}$ consist of one element---the empty function.
Since an increasing function from $[k]$ to $[l]$ is
determined by its range, $\binom{l}{k}$ can be identified with the set of all
$k$ element subsets of $[l]$. Let
${\mathcal F}={\mathcal P}$ be the set of all $\binom{l}{k}$ with
$0<k\leq l$ or $k=l=0$. Declare $\binom{n}{m}\dbullet
\binom{l}{k}$ to be defined precisely when $m=l$, and let
\[
\binom{n}{l}\dbullet \binom{l}{k} = \binom{n}{k}.
\]
It is clear that $\binom{n}{l}\dbullet \binom{l}{k}=
\binom{n}{l}\dowd \binom{l}{k}$. Note, however, that
$\binom{n}{m}\dowd \binom{l}{k}$ is defined if we assume only
$m\geq l$.}
\end{illustration}

The following condition is our Ramsey statement, which is just a formalization
of the statement from the beginning of this subsection:
\begin{enumerate}
\item[{\bf (R)}] given $d>0$, for each $P\in {\mathcal P}$, there is an 
$F\in {\mathcal F}$ such that $F\dbullet P$ is defined, and for
every $d$-coloring of $F\dbullet P$ there is an $f\in F$ such that
$f\dowd P$ is monochromatic.
\end{enumerate}

\begin{illustration} {\rm In the special case of Illustrations~\ref{I:ramfir} and \ref{I:dbul}, condition (R)
says in particular that given $d>0$ and $0<k\leq l$ there exists $m\geq l$ such that
for each $d$-coloring of $\binom{m}{l}\dbullet \binom{l}{k} =
\binom{m}{k}$ there exists $a\in \binom{m}{l}$ such that the set
\[
\{ a \circ x\colon x\in \binom{l}{k}\}
\]
is monochromatic. This is the classical Ramsey theorem.}
\end{illustration}

\subsection{Abstract pigeonhole principle}

We introduce here our pigeonhole principle. The name is purely conventional
as the principle is not a simple abstraction of the well known pigeonhole principle of 
Dirichlet. Rather it is a condition that is easy to check in concrete situations and that implies,
through inductive arguments encoded in Theorem~\ref{T:ram}, the Ramsey condition (R).

We will need an important additional piece of structure. Let $A, X$, and an action $\dowd$
be as above. Let $\partial\colon X\to X$ be a function such that for $a\in A$ and $x\in X$, if
$a\dowd x$ is defined, then $a\dowd\partial x$ is defined and
\begin{equation}\label{E:trun}
\partial (a\dowd x) = a\dowd \partial x.
\end{equation}
Such a function $\partial$ is called a {\bf truncation}. For
$P\subseteq X$, we write
\begin{equation}\label{E:part}
\partial P = \{ \partial x\colon x\in P\}.
\end{equation}

Introduction of the operator $\partial$ equips $X$ with an additional structure and equation 
\eqref{E:trun} states that the action of $A$ on $X$ is implemented by homomorphism of this
structure. In applications to concrete Ramsey theorems, $\partial$ is always a form of 
derivation leading from an object in $X$ to another, less complex object in $X$. In this 
fashion, in proofs, $\partial$ provides a foothold for inductive arguments.

\begin{illustration}\label{I:tru} {\rm We continue the pervious illustrations, in particular, our notation
is as in Illustration~\ref{I:ramfir}. For $x\in X$
with $x\colon [k]\to {\mathbb N}\setminus \{ 0\}$, define
\[
\partial x = x\res [k-1].
\]
(Recall here the convention for the notation $k-1$ adopted in the introduction.) 
It is easy to check that condition \eqref{E:trun} is satisfied.
Note also that, by \eqref{E:part}, $\partial \binom{l}{k}= \binom{l-1}{k-1}$, if $k>1$, and
$\partial \binom{l}{k}= \binom{0}{0}$, if $k\leq 1$. }
\end{illustration}

Let $({\mathcal F}, {\mathcal P}, \dbullet)$ be a pair of families
over $(A, X, \dowd)$ equipped with a truncation $\partial$. We are
ready to formulate our pigeonhole principle. For $P\subseteq X$
and $y\in X$, put
\begin{equation}\label{E:py}
P_{y} = \{ x\in P \colon \partial x = y\}.
\end{equation}
So $P_y$ is the set consisting of those elements of $P$ that
truncate to the same simpler object $y$.
Given $a,b\in A$, we say that $b$ {\bf extends} $a$ if for each $x$ with
$a\dowd x$ defined, we have that $b\dowd x$ is defined and that it is equal to $a\dowd x$.
For $F\in {\mathcal F}$ and $a\in A$, let
\begin{equation}\label{E:ext}
F_a = \{ f\in F\colon f\hbox{ extends }a\}.
\end{equation}

The Ramsey statement (R)
above requires, upon coloring of $F\dbullet P$, stabilizing the
coloring on a copy $f\dowd P$ of $P$ obtained by acting on $P$ by some element
$f$ of $F$. Pigeonhole principle (P) below asks us to
perform the following much easier task. We fix an object $y\in X$,
which can be assumed to be simpler than objects in $P$. We
consider the elements of $P$ that truncate to this fixed $y$, that
is, we consider $P_y$, and require stabilizing the coloring
{\em only} on a copy $f\dowd P_y$ of $P_y$ obtained by acting on $P_y$ by an element $f$
from $F$. The price to pay is that $f$ has to act on
$y$ in a way prescribed by an element $a\in A$ chosen in advance, that is,
$f$ is actually taken from $F_a$ for some $a$ for which $a\dowd y$ is defined.

It is surprising that various concrete pigeonhole principles occurring in the
finite pure Ramsey theory have this form. We illustrate
it below by the classical pigeonhole principle used to prove the
classical Ramsey theorem. In the
following sections, we will give more complex examples involving trees. Paper \cite{Sol}
contains a number of further examples.

The following criterion on
$({\mathcal F}, {\mathcal P}, \dbullet)$ is our pigeonhole
principle:
\begin{enumerate}
\item[{\bf (P)}] given $d>0$, for all $P\in {\mathcal P}$ and $y\in
\partial P$, there are $F\in {\mathcal F}$ and $a\in A$ such
that $F\dbullet P$ is defined, $a\dowd y$ is defined, and for
every $d$-coloring of $F_a\dowd P_{y}$ there is an $f\in F_a$ such that
$f\dowd P_{y}$ is monochromatic.
\end{enumerate}
Note that in the condition above $F_a\dowd P_y$ is defined since $F\dbullet P$ is
assumed to be defined and $F_a\subseteq F$ and $P_y\subseteq P$. Also, of course, 
the condition would not have changed if we required the coloring to be defined on $F\dowd P_y$ or 
even on $F\dowd P$. It is, however, crucial that $f$ be found in $F_a$. 

\begin{illustration}\label{I:P1} {\rm In our special case from the earlier illustrations, a moment of thought
and a picture convince one that condition (P) boils down to the classical
pigeonhole principle. For the sake of practice, however, let us look at it carefully in
detail. We will be helped by Figure~\ref{fig:drawing1}. 

\begin{figure}[htb]
\centering
\def\svgwidth{0.7\columnwidth}
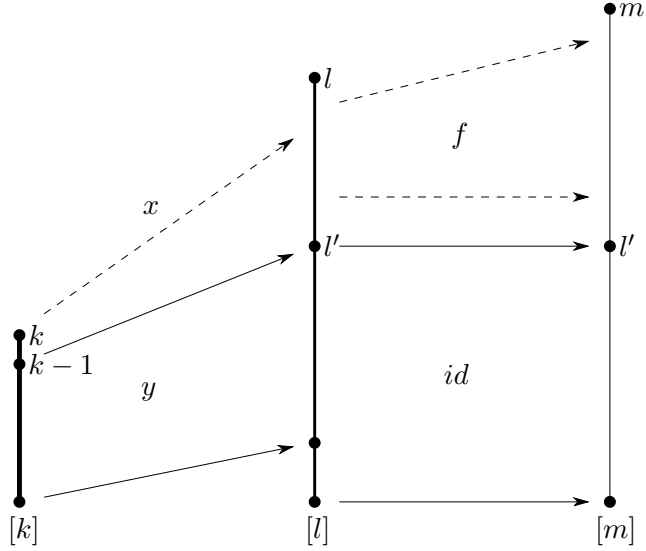
\caption{Condition (P) in Illustration \ref{I:P1}.}
\label{fig:drawing1}
\end{figure}

For notational simplicity, in the argument below, we assume that $k>1$ and leave checking
that the same argument works for $k\leq 1$ to the reader.
We state condition (P) in our special case:

\noindent let $d>0$, $1<k\leq l$ and $y\in \binom{l-1}{k-1}$ be
given; let $l'$ be the maximum of the range of $y$; there exists $m\geq l$ and an increasing 
function $a\colon [l']\to {\mathbb N}\setminus \{ 0\}$ such that for each each $d$-coloring of
\[
\{f\circ x\colon f\in \binom{m}{l},\, f\upharpoonright [l'] = a,\, x\in \binom{l}{k},\,
x\upharpoonright [k-1] = y\}
\]
there exists $f\in \binom{m}{l}$ with $f\upharpoonright [l']
= a$ and such that 
\[
\{ f\circ x\colon x\in \binom{l}{k},\, x\upharpoonright [k-1] = y\}
\]
is monochromatic.

We claim that the condition
above holds with $a$ being the identity function from $[l']$ to
itself. Indeed, with this choice of $a$, the conclusion of the condition reads:

\noindent there exists $m\geq l$ such that for each $d$-coloring of
\[
\{f\circ x\colon f\in \binom{m}{l},\,  f(i)=i\hbox{ for }i\in [l'],\, x\in \binom{l}{k},\,
x\upharpoonright [k-1] = y\}
\]
there exists $f\in \binom{m}{l}$ with $f(i)=i$, for $i\in [l']$,
and with
\[
\{ f\circ x\colon x\in \binom{l}{k},\, x\upharpoonright [k-1] = y\}
\]
monochromatic.

\noindent The elements of the set
\[
\{f\circ x\colon f\in \binom{m}{l},\, f(i)=i\hbox{ for }i\in [l'],\, x\in \binom{l}{k},\,
x\upharpoonright [k-1] = y\}
\]
differ only in the single value $f(x(k))$ and this value comes
from the set $[m]\setminus [l']$. Also $x(k)$ is an arbitrary 
element of $[l]\setminus [l']$. So, in essence, we are
$d$-coloring $[m]\setminus [l']$ and are looking for an increasing
function from $[l]\setminus [l']$ to $[m]\setminus [l']$ whose
values take the same color. This is just the classical pigeonhole
condition, and we can take $m$ to be any number strictly bigger than $l'+d\cdot (l-l'-1)$.}
\end{illustration}

Our goal is to state a theorem that condition (P) implies
condition (R). Achieving this goal, in Theorem~\ref{T:ram}, will
require introducing more structure on $(A, X, \dowd,
\partial)$ and imposing additional conditions on
$({\mathcal F}, {\mathcal P}, \dbullet)$.

\subsection{Additional structure and additional conditions}

Let $A, X$, an action $\dowd$, and a truncation $\partial$ be as
above.

Let again $({\mathcal F}, {\mathcal P}, \dbullet)$ be a pair of
families over $(A,X, \dowd)$. Recall the notion of extension for elements of $A$ defined in
the discussion preceding \eqref{E:ext}. We first state two conditions on
$({\mathcal F}, {\mathcal P}, \dbullet)$ that do not require
introducing any additional structure:
\begin{enumerate}
\item[{\bf (A)}] if $P\in {\mathcal P}$, then $\partial P\in {\mathcal
P}$;

\item[{\bf (B)}] if $F\in {\mathcal F}$, $P\in {\mathcal P}$, and
$F\dbullet\partial P$ is defined, then there is $G\in {\mathcal
F}$ such that $G\dbullet P$ is defined and for each $f\in F$
there is $g\in G$ extending $f$.
\end{enumerate}
Strictly speaking conditions (A) and (B) are
not needed to prove Theorem~\ref{T:ram}; one can dispense with them at
the expense of strengthening condition (P) slightly. (We elaborate on it in Appendix~1.)
However, in some situations, for example, in all the situations in
this note, conditions (A) and (B) hold, and whenever they hold
they do so in an obvious way (and they make strengthening of (P) unnecessary).
Condition (A) simply requires closure of
$\mathcal P$ under truncation. As for condition (B), note that if
$F\dowd P$ is defined, then $F\dowd \partial P$ is defined. The reverse implication is false in general.
Condition (B) gives a substitute for this reverse implication: assuming
something stronger, namely that $F\dbullet
\partial P$ is defined, we can infer that $G\dbullet P$ is defined for a
$G$ that can simulate the action of every element of $F$.

\begin{illustration}
{\rm Recall that we have
\[
{\mathcal F} = {\mathcal P} = \{ \binom{n}{m}\colon 0<m\leq n\hbox{ or }m=n=0\}.
\]
We check conditions (A) and (B). It follows from the remark in Illustration~\ref{I:tru} that
$\mathcal P$ is closed under $\partial$, so (A) holds. To check (B), let $F=\binom{n}{m}$ and $P=\binom{l}{k}$. We assume $k>1$ and
leave the trivial case $k\leq 1$ to the reader. We have
\[
F\dbullet \partial P = \binom{n}{m}\dbullet \binom{l-1}{k-1}
\]
is defined precisely when $m=l-1$, and we can take $G= \binom{n+1}{l}$ to witness the conclusion of (B)
since $\binom{n+1}{l}\dbullet \binom{l}{k}$ is defined and each element of $\binom{n}{l-1}$ is extended by an element of $\binom{n+1}{l}$.
We elaborate on this last point. Note that for each $f\in \binom{n}{l-1}$, that is, for each increasing $f\colon [l-1]\to [n]$, there is
increasing $g\colon [l]\to [n+1]$
with $g\res [l-1] = f$. In this situation, for each $x\in X$ (recall that $X$ is the set of all increasing functions from some $[k]$ to
$\mathbb N$), if $f\dowd x$ is defined, then the image of $x$ is included in $[l-1]$, and so
$g\dowd x$ is defined and obviously
\[
g\dowd x = g\circ x = f\circ x = f\dowd x.
\]
So, in our example, $g$ extending $f$ as an increasing function is equivalent to
$g$ extending $f$ as an element of $A$. A similar coincidence will be present also in the subsequent illustrations.
}
\end{illustration}

To make the partial function $\dowd$ from
$A\times X$ to $X$ into an honest action, we assume that we also have a {\em partial} function from $A\times
A$ to $A$:
\[
(a, b)\to a\cdot b,
\]
such that for $a,b\in A$ and $x\in X$ if $a\dowd (b\dowd x)$ and
$(a\cdot b)\dowd x$ are both defined, then
\begin{equation}\label{E:acct}
a\dowd (b\dowd x) = (a\cdot b)\dowd x.
\end{equation}
The operation $\cdot$ as above will be called {\bf multiplication}. Equation \eqref{E:acct} is the usual equation defining, say,
a group action on a set. As before, for $F,\, G\subseteq A$, we say that
$F\cdot G$ {\em is defined} if $a\cdot b$ is defined for all $a\in F$
and $b\in G$ and we let
\[
F\cdot G = \{ a\cdot b\colon a\in F, b\in G\}.
\]

Now, again as before, in addition to the partial function $\dbullet$ from ${\mathcal F}\times {\mathcal P}$ to $\mathcal P$,
assume that we have a partial function $\bullet$ from ${\mathcal F}\times {\mathcal F}$ to $\mathcal F$
with the property that if $G\bullet F$ is defined, then it is given point-wise, that is, $G\cdot F$
is defined and
\[
G\bullet F = G\cdot F.
\]
We now call $({\mathcal F}, {\mathcal P}, \dbullet, \bullet)$ a {\bf pair of families over}
$(A,X, \dowd, \cdot)$.

We can now state our final condition on $\mathcal F$, $\mathcal P$, $\dbullet$ and $\bullet$:
\begin{enumerate}
\item[{\bf ($*$)}] if $F,G\in {\mathcal F}$, $P\in {\mathcal P}$, and
$F\dbullet (G\dbullet P)$ is defined, then so is $(F\bullet
G)\dbullet P$.
\end{enumerate}
This condition is crucial. It says that $F\dbullet(G\dbullet P)$ is never defined
``by chance;" if it is defined, then the product $F\bullet G$ is defined, as is its action
on $P$. In concrete situations, this condition is guaranteed by a natural calibration
of the domains of the operations $\dbullet$ and $\bullet$.
Note that under the assumptions of ($*$), from \eqref{E:acct}, we have
\[
F\dbullet (G\dbullet P) = (F\bullet G)\dbullet P.
\]
In \cite{Sol}, a pair of families  $({\mathcal F}, {\mathcal P}, \dbullet, \bullet)$ over
$(A,X, \dowd, \cdot)$ fufilling condition ($*$) is called an {\em actoid of sets}.

\begin{illustration}
{\rm Recall again that
\[
{\mathcal F} = {\mathcal P} = \{ \binom{n}{m}\colon 0<m\leq n\hbox{ or }m=n=0\}.
\]
Declare $\binom{n}{m}\bullet\binom{l}{k}$ on $\mathcal F$ to be defined precisely when $m=l$ and let
\[
\binom{n}{l}\bullet\binom{l}{k} = \binom{n}{k}.
\]
So $\bullet$ is equal to $\dbullet$ defined earlier in Illustration~\ref{I:dbul}. It follows that $\bullet$
is given pointwise.

To check ($*$), note that if
\[
\binom{q}{p}\dbullet (\binom{n}{m}\dbullet\binom{l}{k})
\]
is defined, then $m=l$ and $p=n$, but in this situation
\[
(\binom{q}{p}\bullet \binom{n}{m})\dbullet\binom{l}{k}
\]
is defined.
}
\end{illustration}

We require one more piece of structure that, roughly speaking,
measures complexity of objects in $X$. A function $|\cdot | \colon
X\to D$, where $(D,\leq)$ is a linear order, is called a {\bf norm}
if for $x,y\in X$, $|x|\leq |y|$ implies that for all $a\in A$
\begin{equation}\label{E:norm}
a\dowd y \hbox{ defined}\Rightarrow (a\dowd x\hbox{ defined and }
|a\dowd x|\leq |a\dowd y|).
\end{equation}

\begin{illustration}\label{I:ramlas}
{\rm In our special case, $X$ is the set of all increasing injections $x\colon [k]\to {\mathbb N}\setminus \{ 0\}$ 
for $k\in {\mathbb N}$.
Define $|\cdot|\colon X\to {\mathbb N}$, where $\mathbb N$ is taken with its natural linear order, to be
\begin{equation}\notag
|x| =
\begin{cases}
\max {\rm image}(x) = x(k), &\text{ if $k>0$};\\
0, &\text{ if $k=0$}.
\end{cases}
\end{equation}
We check that this definition gives a norm. Let $a\in A$, $a\colon [l]\to {\mathbb N}\setminus \{ 0\}$.
Note that, for $x\in X$, $a\dowd x$ is defined precisely when $|x|\leq l$ and $|a\dowd x| = a(|x|)$, if $|x|>0$, and $|a\dowd x| = 0$, if $|x|=0$.
So given $x_1, x_2\in X$ with $|x_1|\leq |x_2|$,  it is clear that if $a\dowd x_2$ is defined, then so is $a\dowd x_1$ and
\[
|a\dowd x_1|= a(|x_1|)\leq a(|x_2|) = |a\dowd x_2|,\; \hbox{ if }|x_1|>0,
\]
or
\[
|a\dowd x_1|= 0\leq  |a\dowd x_2|,\; \hbox{ if }|x_1|=0.
\] }
\end{illustration}

The additional conditions required to prove our theorem were stated as (A), (B), and ($*$).
The additional structure introduced above is consolidated in the following notion.
\begin{definition}\label{D:norback}
A {\bf normed background} is a pair of sets $A, X$ equipped with a
multiplication $\cdot$ and an action $\dowd$ fulfilling
\eqref{E:acct}, with a truncation $\partial$ fulfilling
\eqref{E:trun}, and with a norm $|\cdot|$ fulfilling
\eqref{E:norm}.
\end{definition}
\noindent With some abuse of notation, a normed background as
above will be denoted by $(A,X)$.

\subsection{The theorem}

Now we can phrase our theorem. To see how it follows from the somewhat more general results of \cite{Sol},
the reader should consult Appendix~1. We write $\partial^tP$, $t\in {\mathbb N}$, for the result of applying
truncation $\partial$ to $P$ $t$ times.

\begin{theorem}\label{T:ram}
Let $({\mathcal F}, {\mathcal P}, \dbullet, \bullet)$ be a pair of families over a normed background fulfilling
conditions (A), (B), and ($*$). Assume that each $P\in {\mathcal P}$ is finite and 
for each $P\in {\mathcal P}$ there
is $t\in {\mathbb N}$ such that $\partial^t P$ consist of one
element. If $({\mathcal F}, {\mathcal P})$ fulfills (P), then it fulfills (R).
\end{theorem}

Note that the theorem above gives the classical Ramsey theorem on
the basis of Illustrations~\ref{I:ramfir}--\ref{I:ramlas}. In them, we checked all the assumptions of Theorem~\ref{T:ram}
except: for $P\in {\mathcal P}$, $P$ is finite and $\partial^tP$ has one element for some $t\in {\mathbb N}$.
Finiteness of $P$ is clear. Note that $\partial^k\binom{l}{k} =\binom{0}{0}$, so this last assumption is also fulfilled. 

\section{Trees and another illustration}\label{S:tree}

\noindent {\bf Trees and embeddings.} We state here basic definitions concerning trees. By a {\bf tree}
we understand a {\em finite}, possibly empty, partial order such that each two elements have a common predecessor and
the set of predecessors of each
element is linearly ordered. So trees for us are {\em finite trees}.
If a tree is non-empty, it has a smallest element, which we call the {\bf root}. Maximal
elements of a tree are called {\bf leaves}. By convention, we regard every node of a tree as one of its
own predecessors and as one of its own successors.

Each tree $T$ carries a binary function $\wedge_T$ that assigns to each
$v,w\in T$ the largest element $v\wedge_Tw$ of $T$ that is a predecessor of
both $v$ and $w$. After Deuber \cite{Deu}, we say that a
function $f\colon S\to T$, for trees $S$ and $T$, is a {\bf morphism} if for all $v,w\in S$,
\[
f(v\wedge_Sw) = f(v)\wedge_T f(w).
\]
So strictly speaking $f$ is a morphism from the functional structure $(S, \wedge_S)$ to the functional structure 
$(T, \wedge_T)$.

For a tree $T$ and $v\in T$, let ${\rm im}_T(v)$ be the set of all {\bf immediate successors} of $v$, and 
we do not regard $v$ as one of them.
Let $T(v)$ be the tree whose elements are all the successors of $v$ (with $v$ among them, of course). Let 
${\rm ht}_T(v)$ be the cardinality of the set of all predecessors of $v$ (including $v$), and let
\[
{\rm ht}(T) = \max \{ {\rm ht}_T(v)\colon v\in T\}.
\]
For a non-empty tree $T$, let
${\rm br}(T)$ be the maximum of cardinalities of ${\rm im}_T(v)$ for $v\in T$.

We will occasionally suppress the subscripts from various pieces of notation introduced above if we 
deem them clear from the context. 

A tree $T$ is called {\bf ordered} if for each $v\in T$ there is a fixed linear order of ${\rm im}(v)$. Such an assignment allows us to define
the lexicographic linear order $\leq_T$ on all the nodes of $T$ by stipulating that $v\leq_T w$ if $v$ is a predecessor of $w$ and, in case $v$ is not a predecessor
of $w$ and $w$ is not a predecessor of $v$, that $v\leq_Tw$
if the predecessor of $v$ in ${\rm im}(v\wedge w)$ is less than or equal to the predecessor of $w$ in 
${\rm im}(v\wedge w)$ in
the given order on ${\rm im}(v\wedge w)$.

The simplest ordered trees are $[n]$ for $n\in {\mathbb N}$ with their natural successor relation and the unique
ordering of the immediate successors of each vertex.

An {\bf embedding $f$ from an ordered tree $S$ to an ordered tree $T$} is an injective
tree morphism such that
\begin{enumerate}
\item[(i)] it is order preserving between $\leq_S$ and $\leq_T$;

\item[(ii)] for each $v\in S$, the set
\[
\{ w\in {\rm im}_T(f(v))\colon w \hbox{ is a predecessor of }f(v')\hbox{ for some }v'\in {\rm im}_S(v)\}
\]
forms an initial segment with respect to $\leq_T$ of ${\rm im}_T(f(v))$.
\end{enumerate}
Note that preservation of order by $f$ is equivalent to saying that
for every $v\in S$ and all $w_1,w_2\in {\rm im}_S(v)$
with $w_1\leq _S w_2$ if $w_1', w_2'$ in ${\rm im}_T(f(v))$ are predecessors of $f(w_1)$ and $f(w_2)$, respectively, then
$w_1'\leq_T w_2'$.
An embedding is {\bf leaf preserving} if each leaf of the domain is mapped to a leaf of the range.
An embedding $f\colon S\to T$ is called {\bf strong} if for $v,w\in S$ with ${\rm ht}(v) = {\rm ht}(w)$ we
have that ${\rm ht}(f(v)) = {\rm ht}(f(w))$. Note that each embedding from $[n]$, $n\in {\mathbb N}$, to an ordered tree
is a strong embedding.

\noindent {\bf Derivations on trees.} There are two natural ways to trim an ordered tree. Let an ordered tree $T$ be given. Put
\begin{equation}\label{E:derst}
T^* = \{ v\in T\colon {\rm ht}(v)<{\rm ht}(T)\},
\end{equation}
that is, $T^*$ is obtained from $T$ by removing all of its highest leaves. Note that $T^*$ with $\leq_T$ restricted to it
is an ordered tree, and that the inclusion from $T^*$ to $T$ is a strong embedding.

Let $x$ be the rightmost with respect to $\leq_T$ leaf of $T$, that is, $x$ is the $\leq_T$-largest element of $T$, and let
\begin{equation}\label{E:derpr}
T' = \{ v\in T\colon T(v) \hbox{ has a leaf different from }x\},
\end{equation}
that is, $T'$ is obtained from $T$ by removing from it a final segment of its rightmost branch.
The tree $T'$ with $\leq_T$ restricted to it forms an ordered tree and the inclusion $T'\subseteq T$ is a leaf preserving
embedding.
If $T\not=\emptyset$, then the set $T\setminus T'$ with the inherited tree structure
can be identified with $[p]$ for some $p\in {\mathbb N}$, $p>0$, with its natural tree order. If $T'\not=\emptyset$, there is a unique node $v_0\in T'$ that
has an immediate successor in $T\setminus T'$. We call $v_0$ the {\bf splitting node of} $T$.

\noindent {\bf Examples of trees and embeddings.} We fix some notation concerning trees. After  
Deuber \cite{Deu}, a non-empty tree
$T$ is called {\bf regular} if for each $v\in T$ that is not a leaf, $|{\rm im}(v)| = {\rm br}(T)$ and for each leaf  
$x\in T$,
${\rm ht}(x) = {\rm ht}(T)$. Of course, each such tree is fully determined by the value of two parameters:
${\rm br}(T)$ and ${\rm ht}(T)$. For $k,n\in {\mathbb N}$, $k>0$, $n>1$, let $T^{k,n}$ be the regular tree of 
height $n$ and with branching number $k$.
By convention, for $k\in {\mathbb N}$, let $T^{k,1}$ have exactly one node and $T^{k,0}$ be equal to the 
empty tree, and for $n\in {\mathbb N}$, $n>1$, let $T^{0,n}$ have exactly one node.
We consider $T^{k,n}$ to be an ordered
tree with some linear order $\leq_{T^{k,n}}$. (All possible orders making $T^{k,n}$ into an ordered 
tree lead to isomorphic ordered
trees.) The tree $T^{1,n}$ can be identified with $[n]$ as an ordered tree.

We fix two natural ways of embedding $T^{k,n}$ into $T^{k, n+1}$.
First, there is a unique embedding $\iota^*$ of $T^{k, n}$ into $T^{k,n+1}$
with
\[
\iota^*({\rm im}_{T^{k,n}}(v)) \subseteq {\rm im}_{T^{k,n+1}}(\iota^*(v)),
\]
for $v\in T^{k,n}$, and with $\iota^*$ mapping the root of $T^{k,n}$ to the root of
$T^{k,n+1}$, if $T^{k,n}$ is non-empty. Note that ${\rm ht}_{T^{k,n}}(v) = {\rm ht}_{T^{k,
n+1}}(\iota^*(v))$ for $v\in T^{k,n}$. We write
\begin{equation}\label{E:inclst}
T^{k,n}\subseteq^* T^{k,n+1}
\end{equation}
to indicate that we consider $T^{k,n}$ identified with its image
under $\iota^*$. There is also a unique embedding $\iota'$ of
$T^{k,n}$ into $T^{k, n+1}$ with
\[
\iota'({\rm im}_{T^{k,n}}(v)) \subseteq {\rm im}_{T^{k,n+1}}(\iota'(v)),
\]
for $v\in T^{k,n}$, and with $\iota'$ mapping the $\leq_{T^{k,n}}$-smallest leaf of $T^{k,n}$ to the
$\leq_{T^{k, n+1}}$-smallest leaf of $T^{k,n+1}$, if $T^{k,n}$ is non-empty. This embedding comes 
from the isomorphism between $T^{k,n}$ and $T^{k,n+1}(v_0)$, where $v_0$ is the 
$\leq_{T^{k,n+1}}$-smallest immediate successor of the root of $T^{k,n+1}$.
Note that the image of the set of all leaves of $T^{k,n}$ under
$\iota'$ is an initial segment with respect to $\leq_{T^{k, n+1}}$
of the set of leaves of $T^{k, n+1}$. We write 
\begin{equation}\label{E:inclpr}
T^{k,n}\subseteq' T^{k,n+1}
\end{equation}
to indicate that we consider $T^{k,n}$ identified with its image under $\iota'$. 

We give one more illustration. Its conclusion will be used in the sequel.
\begin{illustration}\label{I:tree}
{\em We prove the following, possibly folklore, generalization of the classical Ramsey theorem.

\noindent {\em Given $d>0$, $s\in {\mathbb N}$, and a non-empty ordered tree $S$, there is a non-empty 
ordered tree
$T$ with ${\rm br}(T) = {\rm br}(S)$ such that
for each $d$-coloring of all leaf preserving embeddings of $[s]$ to $T$ there exists a leaf preserving embedding $g_0\colon S\to T$
such that
\[
\{ g_0\circ f\colon f\colon [s]\to S\hbox{ a leaf preserving embedding}\}
\]
is monochromatic.}

The proof below consists essentially of stating definitions. All the
checking that needs to be done is routine and would be probably best
left to the reader. However, since this is the first example
involving trees, we will perform all the verifications carefully and
explicitly.

Let $k={\rm br}(S)$. Since there is a leaf preserving embedding from
$S$ to $T^{k,m}$ with $m= {\rm ht}(S)$, we can assume that
$S=T^{k,m}$ for some $m$. For $n\in {\mathbb N}$, set
\[
T^n = T^{k,n}.
\]

\noindent {\bf Normed background.} Let $X$ be the set of all (not
necessarily leaf preserving) embeddings from some $[m]$ to some
$T^n$. Let $A$ consist of all strong embeddings from some $T^m$ to
some $T^n$. (Strong embeddings were defined earlier in this section.
We will need this more restrictive notion of embedding for the
normed background we are defining to work.) For $f\in X$ and
$g,g_1,g_2\in A$ declare that $g\dowd f$ is defined if the image of
$f$ is included with respect to $\subseteq^*$ (as defined in
\eqref{E:inclst}) in the domain of $g$, and similarly declare that
$g_2\cdot g_1$ is defined if the image of $g_1$ is included with
respect to $\subseteq^*$ in the domain of $g_2$, and let
\[
g\dowd f = g\circ f\;\hbox{ and }\; g_2\cdot g_1 = g_2\circ g_1.
\]

Define $\partial^*$ on $X$ be letting for $f\colon [m]\to T^n$,
\[
\partial^* f =  f\res [m-1].
\]
Note, after recalling the derivation \eqref{E:derst}, that $\partial^* f = f\res [m]^*$. Let 
\begin{equation}\label{E:htde}
|f| =
\begin{cases}
{\rm ht}(f(m)),& \text{ if $m>0$};\\
0,& \text{ if $m=0$}.
\end{cases}
\end{equation}
It is checked without any difficulty that $(A,X)$ with the
operations defined above is a normed background. The requirement
that embeddings in $A$ be strong is used in checking that $|\cdot|$
is a norm.

\noindent {\bf A pair of families over $(A,X)$.}
Let $S, T$ be ordered trees. Let
\begin{equation}\label{E:newt}
\binom{T}{S}^s\;\hbox{ and }\;\binomsq{T}{S}^s
\end{equation}
stand for the set of all strong embeddings and strong, leaf
preserving embeddings, respectively, from $S$ to $T$. Since all
embeddings from $[m]$, $m\in {\mathbb N}$, to an ordered tree are
strong, we simplify our notation by setting
\begin{equation}\label{E:emle2}
\binom{T}{m} = \binom{T}{[m]}^s\;\hbox{ and }\;\binomsq{T}{m} = \binomsq{T}{[m]}^s,
\end{equation}
that is, $\binom{T}{m}$ and $\binomsq{T}{m}$ stand for the set of
all embeddings from $[m]$  to $T$ and for the set of all leaf
preserving embeddings from $[m]$ to $T$, respectively.

Let ${\mathcal F}$ consist of sets of the form $\binom{T^n}{T^m}^s$  and
$\binomsq{T^n}{T^m}^s$ with $0<m\leq n$ or $m=n=0$. Declare $\bullet$ to be defined
precisely in the following situations: $\binom{T^n}{T^m}^s\bullet
\binom{T^m}{T^l}^s$ and $\binomsq{T^n}{T^m}^s\bullet
\binomsq{T^m}{T^l}^s$, and define them to be
\[
\binom{T^n}{T^m}^s\bullet \binom{T^m}{T^l}^s=\binom{T^n}{T^l}^s\;\hbox{
and }\; \binomsq{T^n}{T^m}^s\bullet
\binomsq{T^m}{T^l}^s=\binomsq{T^n}{T^l}^s.
\]
Let $\mathcal P$ consist of sets of the form $\binom{T^n}{m}$ and
$\binomsq{T^n}{m}$ with $0<m\leq n$ or $m=n=0$. Declare $\dbullet$ to be defined
precisely in the following situations: $\binom{T^n}{T^m}^s\dbullet
\binom{T^m}{l}$ and $\binomsq{T^n}{T^m}^s\dbullet \binomsq{T^m}{l}$, and let
\[
\binom{T^n}{T^m}^s\dbullet \binom{T^m}{l} =  \binom{T^n}{l}\;\hbox{ and }\;
\binomsq{T^n}{T^m}^s\dbullet \binomsq{T^m}{l} =  \binomsq{T^n}{l}.
\]
It is easy to see that these $\bullet$ and $\dbullet$ when defined are given point-wise. This checking
boils down to showing that each strong embedding from $T^l$ to $T^n$ factors through $T^m$, and the
same for strong, leaf preserving embeddings. Such factorizations are easy to produce. Arguing by induction,  
we see that it is suffices to show their existence for $l<m=l+1\leq n$. Since $l<n$, given a strong (leaf preserving, 
respectively) embedding 
$g\colon T^l\to T^n$, there is $1\leq j\leq n$ such that ${\rm ht}(g(v))\not= j$, for each $v\in T^l$. Fix the largest 
$1\leq i\leq l$ such that ${\rm ht}(g(v))< j$ for all $v\in T^l$ with ${\rm ht}(v)=i$, or let $i=0$ if no such $i$ exists. 
Let $g_1\colon T^l\to T^{l+1}$ be an arbitrary strong (leaf preserving, respectively) embedding such that 
for $v\in T^l$ 
\begin{equation}\notag
{\rm ht}(g_1(v)) = 
\begin{cases}
{\rm ht}(v),&\text{ if ${\rm ht}(v)\leq i$;}\\
{\rm ht}(v)+1,&\text{ if ${\rm ht}(v)\geq i+1$.} 
\end{cases}
\end{equation}  
So there is no element of $T^l$ that gets mapped to a $w\in T^{l+1}$ with ${\rm ht}(w) =i+1$. 
Now it is easy to find a strong (leaf preserving, respectively) embedding $g_2\colon 
T^{l+1} \to T^n$ such that $g= g_2\circ g_1$. (We do it so that ${\rm ht}(g_2(w)) =j$ for all $w\in T^{l+1}$ with 
${\rm ht}(w) = i+1$.)

Note that
\begin{equation}\label{E:trsin}
\partial^* \binom{T^n}{m} = \partial^* \binomsq{T^n}{m} =
\begin{cases}
\binom{T^{n-1}}{m-1},& \text{ if $m>1$;}\\
\binom{T^0}{0},& \text{ if $m\leq 1$.}
\end{cases}
\end{equation}
Using \eqref{E:trsin}, we verify that ${\mathcal F}$ and
${\mathcal P}$ is a pair of families over $(A,X)$ fulfilling
conditions (A) and (B). Condition (A) is clear from \eqref{E:trsin}. We verify condition (B) for $k>1$ in the calculation
below, and leave the trivial case $k\leq 1$ to the reader. Note that
by \eqref{E:trsin}, if
\[
\binom{T^n}{T^m}^s\dbullet \partial^* \binom{T^l}{k} = \binom{T^n}{T^m}^s\dbullet\binom{T^{l-1}}{k-1}
\]
is defined, then $l=m+1$, so $\binom{T^{n+1}}{T^{m+1}}^s\dbullet
\binom{T^l}{k}$ is defined and each $g\in \binom{T^n}{T^m}^s$ is
extended by some $h\in\binom{T^{n+1}}{T^{m+1}}^s$ (that is, for each
$f\in X$ if $g\dowd f$ is defined, then so is $h\dowd f$ and $h\dowd
f  = g\dowd f$); simply view $T^m$ as included in $T^{m+1}$ via
$\subseteq^*$ and take $h\colon T^{m+1}\to T^{n+1}$ to be any strong
embedding with $h\res T^m = g$, that is, $h$ extends $g\colon T^m\to
T^n$ as an embedding. We handle the situation when
\[
\binom{T^n}{T^m}^s\dbullet \partial^* \binomsq{T^l}{k} = \binom{T^n}{T^m}^s\dbullet\binom{T^{l-1}}{k-1}
\]
is defined in the same way, except that in this case $\binomsq{T^{n+1}}{T^{m+1}}^s$ witnesses that (B) holds.

To see condition ($*$), note that if
\[
\binom{T^q}{T^p}^s\dbullet (\binom{T^n}{T^m}^s\dbullet \binom{T^l}{k})
\]
is defined, then $m=l$ and $p=n$, so
\[
(\binom{T^q}{T^p}^s\bullet \binom{T^n}{T^m}^s)\dbullet \binom{T^l}{k})
\]
is defined, as required. We handle the situation when
\[
\binomsq{T^q}{T^p}^s\dbullet (\binomsq{T^n}{T^m}^s\dbullet \binomsq{T^l}{k})
\]
is defined in the same way.

Condition (R) for the above defined pair of families clearly gives the statement from the beginning of 
this illustration. Since for
$\binomsq{T^n}{m}\in {\mathcal P}$ we have
$(\partial^*)^n\binomsq{T^n}{m} = \binom{T^0}{0}$ and $\binom{T^0}{0}$ has exactly one element, by Theorem~\ref{T:ram}, it suffices to
see condition (P).

\noindent {\bf Condition (P).} We carefully check condition (P). Fix an element of $\mathcal P$, which must be of the form
$\binom{T^q}{p}$ or $\binomsq{T^q}{p}$. We consider the first case
first. We assume $p>1$ and leave $p\leq 1$ to the reader. To check (P), recall the pieces of notation set up
in equations \eqref{E:py} and \eqref{E:ext}. Fix $f_0\in \partial^*\binom{T^q}{p} =
\binom{T^{q-1}}{p-1}$.  We need to find an element
$\binom{T^r}{T^q}^s$ of $\mathcal F$ (it suffices, of course, to find
$r$) and $g_0\in A$ so that for each $d$-coloring of
$\binom{T^r}{T^q}^s_{g_0}\dowd \binom{T^q}{p}_{f_0}$ there is $g\in
\binom{T^r}{T^q}^s_{g_0}$ such that $g\dowd  \binom{T^q}{p}_{f_0}$ is
monochromatic.  Note that $\binom{T^r}{T^q}^s\dbullet \binom{T^q}{p}$
is automatically defined.

\begin{figure}[htb]
\centering
\def\svgwidth{\columnwidth}
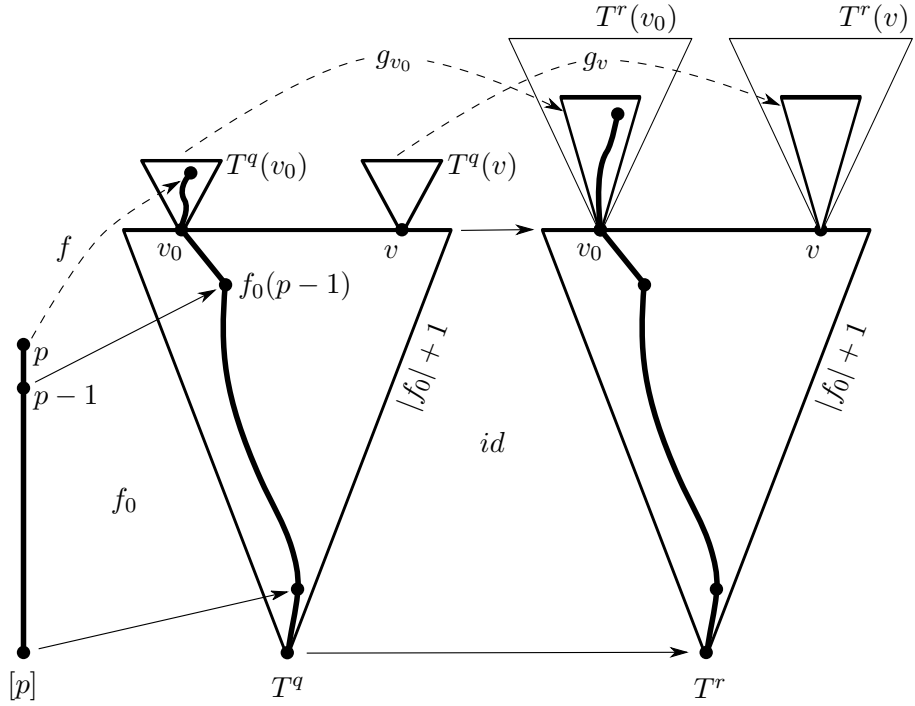
\caption{Condition (P) in Illustration \ref{I:tree}.}
\label{fig:drawing2}
\end{figure}

We claim that $g_0\in A$ equal the identity function
$T^{|f_0|+1}\to T^{|f_0|+1}$ does the job, where $|f_0|$ is defined
by \eqref{E:htde}. Checking (P) boils down to stating precisely what
elements the sets $\binom{T^q}{p}_{f_0}$, $\binom{T^r}{T^q}^s_{g_0}$,
and $g\dowd \binom{T^q}{p}_{f_0}$, for $g\in
\binom{T^r}{T^q}^s_{g_0}$, consist of. Let $v_0$ be the smallest with
respect to $\leq_{T^q}$ element of the set ${\rm im}_{T^q}(f_0(p-1))$, and
keep in mind that we are looking for $r$.

The set $\binom{T^q}{p}_{f_0}$ consists of all $f\in \binom{T^q}{p}$
with $\partial^*f= f_0$. This last condition is equivalent to saying
that $f\res [p-1] = f_0$ and
\begin{equation}\label{E:pig}
f(p)\in T^q(v_0), 
\end{equation}
where \eqref{E:pig} is a consequence of point (ii) in the definition of embedding 
between ordered trees. Each such embedding $f$ is completely determined by the value of $f(p)$.

Fix $r\geq q$, arbitrary for the moment.
Let $g$ be a strong embedding in $\binom{T^r}{T^q}^s_{g_0}$. It is equal to the identity on $T^{|f_0|+1}$ and it is determined by strong
embeddings $g_v$ from $T^q(v)$ to $T^r(v)$,
where $v$ varies over the nodes of $T^q$ with ${\rm ht}(v) = |f_0|+1$. Now, elements of
$g\dowd \binom{T^q}{p}_{f_0}$ are embeddings $g\circ f\colon [p]\to T^r$ with $f$ for which \eqref{E:pig} holds.
Each such embedding is completely determined by the value
\[
(g\circ f)(p) = g_{v_0}(f(p))\in T^r(v_0).
\]
Therefore, solving the problem of fixing the color on $g\dowd
\binom{T^q}{p}_{f_0}$ amounts to the following: $d$-color $T^r(v_0)$
(this is where the values of $(g\circ f)(p)$ are coming from), then
find a strong embedding (this is $g_{v_0}$) of $T^q(v_0)$ (this is where
the values $f(p)$ are located) to $T^r(v_0)$  so that the image of
$T^q(v_0)$ is monochromatic. This can be arranged using a form of the Halpern--L{\"a}uchli theorem 
((HL2) with $t=1$ from Appendix~2) by taking $r$ large enough since $T^q(v_0)$
and $T^r(v_0)$ are isomorphic to $T^m$ and $T^n$, where
$m=q-(|f_0|+1)$ and $n=r-(|f_0|+1)$, respectively. For $v\not= v_0$, after identifying
$T^q(v)$ with $T^q(v_0)$ and $T^r(v)$ with $T^r(v_0)$ via the unique isomorphisms, we let
$g_v$ be equal to $g_{v_0}$. Note that so defined $g$ is strong.

The case $P=\binomsq{T^q}{p}$ is handled analogously with the exception that for $F$ one takes $\binomsq{T^r}{T^q}^s$ for large enough
$r$ and one uses another form of the Halpern--L{\"a}uchli theorem ((HL1) from Appendix~2). We leave it to the reader 
to re-check the details.}
\end{illustration}

\section{A Ramsey theorem for finite trees}\label{S:JD}

We prove the following theorem that extends the results of Deuber \cite{Deu} and of Jasi{\'n}ski
\cite{Jas}. Our proof differs from the arguments of these two papers.

\begin{proposition}\label{P:gen} For non-empty ordered trees $S,T$ and $d>0$, there exists a non-empty
ordered tree $V$ with ${\rm br}(V)={\rm br}(T)$ such that for each
$d$-coloring of all leaf preserving embeddings from $S$ to $V$ there
is a leaf preserving embedding $g_0\colon T\to V$ such that
\[
\{ g_0\circ f\colon f\colon S\to T\hbox{ a leaf preserving
embedding}\}
\]
is monochromatic.
\end{proposition}

As a direct consequence of the above result, one gets its version for embeddings that are not
necessarily leaf preserving by the following argument. Given ordered trees $S,T$, let $S_+, T_+$ be the trees
obtained from $S$ and $T$ by adding one node on top of each leaf of $S$ and $T$, respectively.
Apply now the above statement to $S_+, T_+$ obtaining $V$. Let $V_-$ be gotten from $V$ by deleting
from it all of its leaves. It is easy to check that $V_-$ works by using the obvious observation that embeddings
from $S$ to $T$, from $S$ to $V_-$, and from $T$ to $V_-$  are precisely restrictions of leaf preserving embeddings
from $S_+$ to $T_+$, from $S_+$ to $V$, and from $T_+$ to $V$, respectively.

Deuber's theorem \cite{Deu} is the above result for embeddings that are not necessarily leaf preserving and under the additional
assumptions that ${\rm br}(S) = {\rm br}(T)$ and that $S$ is regular as defined in Section~\ref{S:tree}.
Jasi{\'n}ski's theorem is originally \cite{Jas} stated for boron structures as defined in \cite{Cam}, but can easily be
rephrased in terms of trees, and then it becomes equivalent to the above result
with the additional assumptions that ${\rm br}(S) = {\rm br}(T)$ and that for
each $v\in S$ that is not a leaf $|{\rm im}(v)| = 2$.

We show now how to derive Proposition~\ref{P:gen} from Theorem~\ref{T:ram}.

Let $k={\rm br}(T)$ and set $T^n= T^{k,n}$. Note that it is enough to prove the theorem for $T$ equal to some $T^n$ since
every ordered tree $T$ with ${\rm br}(T)=k$ embeds leaf-preservingly into some $T^n$.

We define an analogue of the set of natural numbers for the present Ramsey situation. We view $T^n$ as an ordered subtree $T^{n+1}$ via the inclusion
$\subseteq'$ defined by \eqref{E:inclpr}. This convention gives an increasing sequence $(T^n)_{n\in {\mathbb N}}$ of ordered trees. Let the
direct limit (that is, the union, if $T^n$ is identified with its image in $T^{n+1}$) of this sequence
be denoted by $T^\infty$. Observe that $T^\infty$ carries a linear order induced from the linear orders $\leq_{T^n}$ on the $T^n$-s.
We denote this linear order by $\leq^\infty$. Each element $v$ of $T^\infty$ belongs to some $T^n$. We call $v$
a {\bf leaf} if $v$ is a leaf in some, or equivalently all, $T^n$ to which it belongs. For an ordered tree $S$, each function $f\colon S\to T^\infty$
has its range included in some $T^n$. We call $f$ a {\bf leaf preserving embedding} if $f$ is a leaf
preserving embedding to some, or equivalently all, $T^n$ in which the image of $f$ is included.
Further, $g\colon D\to T^\infty$ for a subset $D$ of $T^\infty$ is called a {\bf leaf preserving embedding} if the restriction
of $g$ to each $D\cap T^n$, $n\in {\mathbb N}$, is a leaf preserving embedding, where $D\cap T^n$ is taken with the tree order
inherited from $T^n$. For a leaf $x\in T^\infty$, let
\[
T_x = \{ v\in T^\infty\colon v\leq^\infty x\}.
\]
Note that $T_x$ is an infinite set.

\noindent {\bf Normed background.} Let $Y$ consist of all leaf preserving embeddings $f\colon S\to T^\infty$, where $S$ is an ordered tree.
Let $B$ consist of the empty function and of all leaf preserving embeddings $g\colon T_x\to T^\infty$, where $x$ is a leaf of $T^\infty$. It
is easy to see that for such a $g\colon T_x\to T^\infty$, we have $g(T_x)\subseteq T_{g(x)}$. As always, for $f\in Y$ and
$g\in B$, let $g\dowd f$ to be defined precisely when the image of $f$ is included in the domain of $g$ and let
\[
g\dowd f = g\circ f.
\]
Similarly for $g_1, g_2\in B$, define $g_2\cdot g_1$ to be defined precisely when the image of $g_1$ is contained in the
domain of $g_2$ and let
\[
g_2\cdot g_1 = g_2\circ g_1.
\]
We define a truncation using the branch cutting derivation on trees given by \eqref{E:derpr}.
For $f\in Y$ with $f\colon S\to T^\infty$ define
\[
\partial' f = f\res S',
\]
where $S'$ is given by \eqref{E:derpr}.
We define a norm $|\cdot |\colon Y\to T^\infty\cup \{ -\infty\}$, where $T^\infty$ is considered as a linear order with $\leq^\infty$ and $-\infty$ is an element that is less
than all the elements
of $T^\infty$, by letting for $f\in Y$
with $f\colon S\to T^\infty$
\begin{equation}\notag
|f|=
\begin{cases}
\max {\rm image}(f),& \text{ if $S\not=\emptyset$};\\
-\infty,& \text{ if $S=\emptyset$.}
\end{cases}
\end{equation}
Observe that if $S\not=\emptyset$, then $|f|$ is the $\leq^\infty$-minimal leaf $x\in T^\infty$ such that ${\rm image}(f)\subseteq T_x$.
It is easy to check that with so defined operations, $(B, Y)$ becomes a normed background.

\noindent {\bf A pair of families over $(B,Y)$.}
For $n\in {\mathbb N}$, let $x_n\in T^\infty$ be the rightmost leaf of $T^n$ and let $v_n\in T^\infty$ be the root of $T^n$.
Note that
\[
T_{x_n} = T^n \cup \{ v_{n+k}\colon k\in{\mathbb N},\, k>0\}.
\]
Define for $0< m\leq n$
\begin{equation}\notag
\begin{split}
\binomsq{T^n}{T^m}^\infty = \{ g\in B\colon &g\colon T_{x_m}\to T_{x_n},\, g(T^m)\subseteq T^n,\, \hbox{ and }\\
&g(v_{m+k}) =v_{n+k} \hbox{ for all } k\in {\mathbb N}, k>0\}.
\end{split}
\end{equation}
Additionally, let  $\binomsq{T^0}{T^0}^\infty$ consist of the empty function. Observe that the function
\[
\binomsq{T^n}{T^m}^\infty \ni g\to g \res T^m
\]
is a bijection from $\binomsq{T^n}{T^m}^\infty$ to all leaf preserving
embeddings from $T^m$ to $T^n$.

Let ${\mathcal G}$ consist be the family of all subsets of $B$ of the form
$\binomsq{T^n}{T^m}^\infty$, where $n,m\in {\mathbb N}$ and $0< m\leq n$ or $m=n=0$.
Let ${\mathcal Q}$ be the family of all non-empty finite sets $Q\subseteq Y$ of the following form: there is an ordered tree
$S$ such that $Q$ consists of some leaf preserving
embeddings from $S$ to $T^\infty$. In such a situation, we say that $Q$ {\bf is based on} $S$.
As usual, declare $\binomsq{T^n}{T^m}^\infty\bullet \binomsq{T^l}{T^k}^\infty$
to be defined precisely when $m=l$ and let
\[
\binomsq{T^n}{T^l}^\infty\bullet \binomsq{T^l}{T^k}^\infty = \binomsq{T^n}{T^k}^\infty.
\]
Declare $\binomsq{T^n}{T^m}^\infty\dbullet Q$
to be defined precisely when $m$ is the smallest natural number with the property that the
images of all elements of $Q$ are included in $T^m$, and let
\[
\binomsq{T^n}{T^m}^\infty\dbullet Q = \binomsq{T^n}{T^m}^\infty\dowd Q.
\]

We leave to the reader the easy check that $({\mathcal F}, {\mathcal Q},\dbullet,\bullet)$ is
a pair of families over $(B,Y)$, that is, the operations $\dbullet$ and $\bullet$ are given pointwise. The pair of families
fulfills conditions (A), (B), and ($*$). Condition (A) is clear.
To see condition (B), assume that $\binomsq{T^n}{T^m}^\infty\dbullet \partial Q$ is defined, that is,
$m$ is smallest such that the image of all elements of $\partial Q$ is included in $T^m$. Let
$l\in {\mathbb N}$ be smallest such that the image of each element of $Q$ is included in $T^{m+l}$.
Then  $\binomsq{T^{n+l}}{T^{m+l}}^\infty$ witnesses that (B) holds since $\binomsq{T^{n+l}}{T^{m+l}}^\infty\dbullet Q$ is defined
and, as is easy to check, each leaf preserving embedding from $\binomsq{T^n}{T^m}^\infty$ extends (as a function) to a leaf preserving embedding
from $\binomsq{T^{n+l}} {T^{m+l}}^\infty$.
Condition ($*$) follows immediately from an easy observation that if $m$ is the smallest natural number such that the image of each
function in $Q$ is included in $T^m$, then $n$ is the smallest natural number with each function in $\binomsq{T^n}{T^m}^\infty\dowd Q$
having its image included in $T^n$.

Note that condition (R) in this case is the theorem we are proving. Observe also that if $Q\in {\mathcal Q}$ is based on $S$, then
$\partial Q$ is based on $S'$, and $S'$ has one leaf fewer than $S$ if $S\not=\emptyset$. Thus,
$\partial^t Q$ has exactly one element (the empty function) for $t$ equal to the number of leaves in $S$. It follows that to get (R) it remains to check condition (P).

\noindent{\bf Condition (P).} Let $Q\in {\mathcal Q}$ be based on
$S$ and let $q\in {\mathbb N}$ be smallest such that all elements of
$Q$ have ranges included in $T^q$. The set $\partial Q$ is based on
$S'$. We assume $S'$ is not the empty tree. (The case $S'=\emptyset$
is easier, and we ask the reader to handle it after reading the
current argument.) Let $u_0\in S'$ be the splitting node of $S$, and
identify $S\setminus S'$ with $[p]$ for some non-zero $p\in {\mathbb
N}$. (Recall here the discussion following \eqref{E:derpr}.) Fix
$f_0\in
\partial' Q$. Then $f_0\colon S' \to T^q$, $f_0\in Y$. Let
\begin{equation}\label{E:vee}
v_0= f_0(u_0)\in T^q.
\end{equation}

To check (P), we need to find $r\in {\mathbb N}$ and $g_0\in B$ such
that for each $d$-coloring of $\binomsq{T^r}{T^q}^\infty_{g_0}\dowd
Q_{f_0}$ there is $g\in  \binomsq{T^r}{T^q}^\infty_{g_0}$ with $g\dowd
Q_{f_0}$ monochromatic. We will show that large enough $r$ works.
Fix $r\geq q$. Now, we define $g_0$. Find the $\leq^\infty$-smallest
leaf $x$ in $T^q$ such that the image of $f_0$ is included in
$T_x$. Note that $v_0$ is a predecessor of $x$.

First we define $g_0\colon T_x\to T^\infty$. Note that
\[
T_x = (T_x\cap T^q)\cup \{ v_{q+k}\colon k\in {\mathbb N},\, k>0\}.
\]
For the moment, we view $T^q$ as
a subset of $T^r$ in the sense $T^q\subseteq^* T^r$, as defined by
\eqref{E:inclst}, and we let $g_0$ be the identity on the elements of $T_x\cap T^q$ that are not leaves.
Let $g_0$ map leaves of $T_x\cap T^q$ to leaves of $T^r$ in such a way that $g_0$ on $T_x\cap T^q$ is a leaf preserving
embedding to $T^r$. Let $g_0(v_{q+k}) = v_{r+k}$ for $k\in {\mathbb N}$, $k>0$. It is clear that $g_0\in B$.

Consider the set $E$ of all $w\in T^q$ such that $w$ is an immediate
successor of a predecessor of $x$ and $x<^\infty w$. The set
$T^q\setminus T^q_x$ is partitioned into trees $T^q(w)$ with $w\in
E$. Therefore, each $g\in \binomsq{T^r}{T^q}^\infty_{g_0}$ is equal to
$g_0$ in $T_x$ and is completely determined by leaf preserving
embeddings
\[
g_{w}=g\res T^q(w)\colon T^q(w) \to T^r(w),\; w\in E.
\]

Note that $v_0$ given by \eqref{E:vee} has an immediate successor in
$E$. Let $w_0$ be the $\leq^\infty$-smallest among them. For each
$f\in Q_{f_0}$, $f\res S'$ is equal to $f_0$ whose image is included
in $T_x$, while the image of $f\res (S\setminus S')$ is included
in $T^q(w_0)$. So each element of $g\dowd Q_{f_0}$ being of the form
$g\circ f\colon S\to T^r$ is completely determined by
\[
g_{w_0}\circ (f\res (S\setminus S'))\colon S\setminus S' \to T^r(w_0).
\]

Note that the identification of $S\setminus S'$ with $[p]$ makes $f\res (S\setminus S')$ into a leaf
preserving embedding from $[p]$ to $T^q(w_0)$.
Thus, fixing the color on $g\dowd Q_{f_0}$ amounts to the following (with notation as in \eqref{E:emle2}):
$d$-color $\binomsq{T^r(w_0)}{p}$, find a leaf preserving embedding
$g_{w_0}\colon T^q(w_0) \to T^r(w_0)$ so that
$g_{w_0}\dowd \binomsq{T^q(w_0)}{p}$
is monochromatic. This can be achieved from Illustration~\ref{I:tree} by taking $r$ large enough as
$T^r(w_0)$ and $T^q(w_0)$ are isomorphic to $T^n$ and $T^m$, where
$n= r-{\rm ht}(w_0)$ and $m = q-{\rm ht}(w_0)$. We can let $g_{w}\colon T^q(w) \to T^r(w)$ be
arbitrary leaf preserving embeddings for $w\in E$, $w\not= w_0$.

\section{Milliken's theorem in exercises}\label{S:M}

We prove in this section the following result due to Milliken \cite{Mil}. The reader may consult
\cite{Soc} for another purely finitary proof of Milliken's theorem.

\begin{proposition}\label{P:mille} Let $S$ and $T$ be ordered trees. Assume that all leaves in $T$ have the same height.
For $d>0$, there exists an
ordered tree $V$ with ${\rm br}(V)={\rm br}(T)$ such that for each
$d$-coloring of all strong, leaf preserving embeddings from $S$ to
$V$ there is a strong, leaf preserving embedding $g_0\colon T\to V$
such that
\[
\{ g_0\circ f\colon f\colon S\to T\hbox{ a strong, leaf preserving
embedding}\}
\]
is monochromatic.
\end{proposition}

The proof of Proposition~\ref{P:mille} that we will give yields also the statement obtained from Proposition~\ref{P:mille} by replacing strong,
leaf preserving embeddings by strong embeddings in all places. This statement can also be obtained from
Proposition~\ref{P:mille} by a proof that is identical to the argument following Proposition~\ref{P:gen}.
It suffices to notice that, with the notation as in that argument,  strong embeddings
from $S$ to $T$, from $S$ to $V_-$, and from $T$ to $V_-$  are precisely restrictions of strong, leaf preserving
embeddings from $S_+$ to $T_+$, from $S_+$ to $V$, and from $T_+$ to $V$, respectively.

The proof of Proposition~\ref{P:mille} is a somewhat more sophisticated version of the argument in Illustration~\ref{I:tree}.
Let $k={\rm br}(T)$. As before set $T^n = T^{k,n}$. View $T^n$ as a subtree of $T^{n+1}$ via the inclusion
$\subseteq^*$ defined in \eqref{E:inclst}. This inclusion is a strong embedding. This way we obtain an increasing sequence
$(T^n)_{n\in {\mathbb N}}$ of ordered trees. Let $T_\infty$ be the union (direct limit) of this sequence.
The range of each function $f\colon S\to T_\infty$ on an ordered tree $S$
is included in some $T^n$. We call $f$ a {\bf strong embedding} if $f$ is a strong embedding as a function from $S$
to $T^n$ for some, or equivalently all, $T^n$ in which the image of $f$ is included. For $v\in T_\infty$,
let ${\rm ht}(v)$ be equal to ${\rm ht}_{T^n}(v)$ for some, or equivalently, all $T^n$ with $v\in T^n$.

\noindent {\bf Normed background.} Let $Z$ consist of all strong embeddings $f\colon S\to T_\infty$, where $S$ is an ordered tree.
Let $C$ consist of all strong embeddings $g\colon T^m\to T^n$, for some $m\leq n$.
For $f\in Z$ and $g\in C$,
let $g\dowd f$ be defined precisely when the image of $f$ is included in the domain of $g$ and let
\[
g\dowd f = g\circ f.
\]
Similarly for $g_1, g_2\in C$, let $g_2\cdot g_1$ be defined precisely when the image of $g_1$ is contained in the
domain of $g_2$, and let
\[
g_2\cdot g_1 = g_2\circ g_1.
\]
For $f\in Z$ with $f\colon S\to T^n$ define
\[
\partial^* f = f\res S^*.
\]
Define a norm $|\cdot |\colon Z\to {\mathbb N}$, by letting for $f\in Y$
with $f\colon S\to T_\infty$
\[
|f|= \max_{v\in S} {\rm ht}(f(v)).
\]

\noindent {\bf Exercise.} Check that $(C, Z)$ is a normed background.

\noindent {\bf A pair of families over $(C,Z)$.}
The pair of families described below extends the one described in Illustration~\ref{I:tree}.
Recall the sets $\binom{T}{S}^s$ and $\binomsq{T}{S}^s$ defined in equation \eqref{E:newt}.
Let ${\mathcal H}$ consist of all $\binom{T^n}{T^m}^s$ and $\binomsq{T^n}{T^m}^s$ where $m,n\in {\mathbb N}$ and
$0<m\leq n$ or $m=n=0$. Let ${\mathcal R}$ consist of all non-empty sets of the
form $\binom{T^n}{S}^s$ and $\binomsq{T^n}{S}^s$, where $S$ is an ordered tree. Declare
$\binom{T^n}{T^m}^s\bullet \binom{T^l}{T^k}^s$ and $\binomsq{T^n}{T^m}^s\bullet \binomsq{T^l}{T^k}^s$ to be defined
precisely when $m=l$, and let
\[
\binom{T^n}{T^l}^s\bullet \binom{T^l}{T^k}^s = \binom{T^n}{T^k}^s \;\hbox{ and }\;\binomsq{T^n}{T^l}^s\bullet \binomsq{T^l}{T^k}^s = \binomsq{T^n}{T^k}^s.
\]
Similarly, declare $\binom{T^n}{T^m}^s\dbullet \binom{T^l}{S}^s$ and  $\binomsq{T^n}{T^m}^s\dbullet \binomsq{T^l}{S}^s$
to be defined precisely when $m=l$, and let
\[
\binom{T^n}{T^l}^s\dbullet \binom{T^l}{S}^s= \binom{T^n}{S}^s\;\hbox{ and }\;
\binomsq{T^n}{T^l}^s\dbullet \binomsq{T^l}{S}^s= \binomsq{T^n}{S}^s.
\]
The operations $\bullet$ and $\dbullet$ are undefined in situations
not specified above.

\noindent {\bf Exercise.} Check that $({\mathcal H}, {\mathcal
R},\dbullet,\bullet)$ is a pair of families over $(C,Z)$ fulfilling
conditions (A), (B), and ($*$). (Hint. This is almost identical to the argument in Illustration~\ref{I:tree}.)

\noindent {\bf Exercise.} Note that it suffices to prove Proposition~\ref{P:mille} for $T$ of the form $T^n$ (this is where the assumption
that all leaves in $T$ have the same height is used) and check that condition (R) for
$({\mathcal H}, {\mathcal R}, \dbullet, \bullet)$ implies Proposition~\ref{P:mille}
(as well as the statement obtained from Proposition~\ref{P:mille} by replacing strong, leaf preserving embeddings by strong embeddings).

\noindent {\bf Exercise.} Check condition (P) for $({\mathcal H},
{\mathcal R},\bullet,\dbullet)$. (Hint. This follows from the Halpern--L\"{a}uchli theorem for strong subtrees
(HL1) and (HL2) from Appendix~2 and is similar to the argument for (P) in Illustration~\ref{I:tree}.)

\section{Appendix 1: conditions (A) and (B) removed and the final word on normed backgrounds}\label{S:AB}

{\bf 1.} The following criterion (${\rm P}+$)
is the strengthening of condition (P) allowing us to get rid of
conditions (A) and (B). It is obtained from (P) by replacing all occurrences of $P$, except the one in $F\dbullet P$,
by $\partial^tP$ for a fixed but arbitrary $t\in {\mathbb N}$.
\begin{enumerate}
\item[{\bf (P$+$)}] given $d>0$ and $t$, for all $P\in {\mathcal P}$ and $x\in
\partial^{t+1} P$, there are $F\in {\mathcal F}$ and $a\in A$ such
that $F\dbullet P$ is defined, $a\dowd x$ is defined, and for
every $d$-coloring of $F_a\dowd (\partial^tP)_{x}$ there is $f\in F_a$ such that
$f\dowd (\partial^tP)_{x}$ is monochromatic.
\end{enumerate}

The following result is \cite[Corollary 4.4]{Sol}.

\begin{theorem}\label{T:ramful}
Let $({\mathcal F}, {\mathcal S}, \dbullet, \bullet)$ be a pair of families with ($*$) over a normed background.
Assume that each $P\in {\mathcal P}$ is finite and for each $P\in {\mathcal P}$ there
is $t\in {\mathbb N}$ such that $\partial^t P$ consist of one
element. If $({\mathcal F}, {\mathcal P})$ fulfills (P+), then it
fulfills (R).
\end{theorem}

To see that Theorem~\ref{T:ram} is a consequence Theorem~\ref{T:ramful}, we note that (P) in the presence of (A) and (B)
implies (${\rm P}+$).
To see this implication, we proceed by induction on $t$. Condition (${\rm P}+$) for $t=0$ is just (P).
Assuming that (${\rm P}+$) holds for $t$, we prove it for $t+1$. Let $P\in {\mathcal P}$ and $x\in \partial^{t+2}P$.
By condition (A), $\partial P\in {\mathcal P}$. So condition (${\rm P}+$) for $t$ applied to $\partial P$ and $x$ gives
$F\in {\mathcal F}$ and $a\in A$ such
that $F\dbullet \partial P$ is defined, $a\dowd x$ is defined, and for
every $d$-coloring of $F_a\dowd (\partial^{t+1}P)_{x}$ there is $f\in F_a$ such that
$f\dowd (\partial^{t+1}P)_{x}$ is monochromatic. Now condition (B) gives $G\in {\mathcal F}$ such that $G\dbullet P$
is defined and such that each element of $F$ is extended by an element of $G$. It follows that each element of
$F_a$ is extended by an element of $G_a$. Now it is clear that $G$ and $a$ witness that (${\rm P}+$) holds for
$t+1$.

{\bf 2.} The main algebraic structures in the paper are
normed backgrounds. We list below conditions that are more symmetric
than those defining normed backgrounds. As
indicated by Lemma~\ref{L:back}, they give a notion that is in essence equivalent to normed background. All
the normed backgrounds in the present paper and in \cite{Sol}
fulfill the conditions below.

Let $(A, X, \cdot, \dowd, \partial, |\cdot|)$ be such that $\cdot$ is a partial function from $A\times A$ to $A$,
$\dowd$ is a partial function from $A\times X$ to $X$, $\partial$ is a function from $X$ to $X$ and $|\cdot |$ is
a function from $X$ to a set with a linear order $\leq$. Assume the following axioms hold for all $a,b\in A$ 
and $x,y\in X$:
\begin{enumerate}
\item[(i)] if $a\dowd (b\dowd x)$ and $(a\cdot b) \dowd x$ are defined, then
$a\dowd (b\dowd x) = (a\cdot b) \dowd x$;

\item[(ii)] if $a\dowd x$ and $a\dowd \partial x$ are defined, then
$\partial(a\dowd x) = a\dowd \partial x$;

\item[(iii)] $|\partial x|\leq |x|$;

\item[(iv)] if $|x|\leq |y|$ and  $a\dowd x$ and $a\dowd y$ are defined,
then $|a\dowd x|\leq |a\dowd y|$;

\item[(v)] if $|x|\leq |y|$ and $a\dowd y$ is defined, then so is $a\dowd x$.
\end{enumerate}

The following result is \cite[Lemma~4.5]{Sol}.

\begin{lemma}\label{L:back}
\begin{enumerate}
\item[(a)] Assume $(A, X, \cdot, \dowd, \partial, |\cdot|)$ fulfills conditions (i)--(v) above, then $(A,X)$ with $\cdot$, $\dowd$, $\partial$
and $|\cdot|$ is a normed background

\item[(b)] If $(A,X)$ with $\cdot$, $\dowd$, $\partial$ and $|\cdot|$ is a normed background, then there is a function $|\cdot|_1$ on $X$
such that $(A, X, \cdot, \dowd, \partial, |\cdot|_1)$ fulfills conditions (i)--(v) above.
\end{enumerate}
\end{lemma}

\section{Appendix 2: The Halpern--L\"{a}uchli theorem for strong subtrees
as a restatement of the Hales--Jewett theorem}

We point out here that the Halpern--L\"{a}uchli theorem for strong subtrees (there are other, more difficult,  versions) and 
the Hales--Jewett theorem are identical statements phrased in different languages. The 
importance of this translation for the presentation here comes from the fact that the Hales--Jewett theorem 
is shown in \cite{Sol} to be one of the results that follow from the abstract approach to Ramsey theory. 
So when using the Halpern--L\"{a}uchli theorem in the present paper we stay within this approach. 
Justin Moore remarks that
equivalence of these two statements (that is, of the Hales--Jewett and the Halpern--L\"{a}uchli theorems) 
has been known for some time.

We set up a dictionary for translating the Hales--Jewett theorem to the Halpern--L\"{a}uchli theorem.
Let $S$ and $T$ be ordered trees. Let $f\colon {\rm leaves}(S) \to {\rm leaves}(T)$ 
be strictly increasing with respect the orders $\leq_S$ and $\leq_T$ (restricted to 
the leaves),
and be such that for each $v\in S$ there is $w\in T$ such that for
any two leaves $x,y$ of $S$ with $v=x\wedge y$ we have
$w= f(x)\wedge f(y)$. Then there is a unique leaf preserving
embedding from $S$ to $T$ whose restriction to ${\rm leaves}(S)$ is
equal to $f$. We, therefore, refer to such an $f$ itself as a {\bf leaf preserving
embedding}. If in the above definition ${\rm ht}(w)$ depends only
on ${\rm ht}(v)$, then the induced embedding is strong and again
we call $f$ {\bf strong}. A sequence $f_1, \dots, f_r\colon {\rm leaves}(S)\to {\rm leaves}(T)$ of
strong embeddings is called a {\bf strong sequence} if for $x,y\in {\rm leaves}(S)$ and 
$1\leq i,j\leq r$
\[
{\rm ht}(f_i(x)\wedge f_i(y)) = {\rm ht}(f_j(x)\wedge f_j(y)).
\]

Fix a linearly ordered finite set $A$ that is disjoint from $\mathbb
N$. For $n\in {\mathbb N}$, we consider ordered trees
\[
A^{\leq n} = \{ v\colon [l]\to A\colon l\leq n\},
\]
where the tree relation is equal to the extension relation and the
order relation is the one coming from the linear order on $A$. So $A^{\leq n}$ is a version of
the trees $T^{k,n}$ defined in Section~\ref{S:tree}, where $k=|A|$. Note
that the set of leaves of this tree is equal to the set $A^n$ of all
functions from $[n]$ to $A$.

For any function $v\colon [l]\to A$, let $v'\colon A\cup [l]\to A$
be equal to the identity function on $A$ and to $v$ on $[l]$.
Assume we have a function $w\colon [n]\to A\cup [m]$ such that
\begin{enumerate}
\item[(i)] $[m]$ is included in the image of $w$;

\item[(ii)] $w([l])\cap [m]$ is an initial segment of $[m]$,  for each $l\leq n$.
\end{enumerate}
Such $w$ gives rise to a strong embedding $g_w\colon A^m\to A^n$ (recall that $A^m$ and $A^n$
are the sets of leaves of $A^{\leq m}$ and $A^{\leq n}$, respectively)
defined by
\[
g_w(x) = x'\circ w,
\]
for $x\in A^m$. It is easy to check, using property (ii) of $w$,
that $g_w$ preserves the lexicographic order. Property (i) ensures
that $g_w$ is injective. Further note that for $x,y\in A^m$, if
$x\wedge y = v_0$ with ${\rm ht}(v_0) = i_0$, then
$g_w(x)\wedge g_w(y) = v_1$ with ${\rm ht}(v_1) = i_1$,
where
\begin{equation}\label{E:isvs}
i_1 = \max\{ i\colon w([i])\cap [m]\subseteq  [i_0]\}\;\hbox{ and } \;v_1(i) = v_0'( w(i)), \hbox{ for }i\in [i_1].
\end{equation}
Note that $i_1$ depends only on $i_0$. Thus, $g_w$ is indeed a strong embedding.

Now assume that for $r\in {\mathbb N}$, we have $w\colon [n]\to
A^r\cup [m]$ with properties (i) and (ii) above. Such a $w$ gives
rise to $r$ functions $w_i = \pi_i\circ w$, where $\pi_i\colon A^r\cup
[m]\to A\cup [m]$ is the $i$-th projection on $A^r$ and the identity
on $[m]$, also fulfilling conditions (i) and (ii). We therefore get a sequence of
strong embeddings $g_w^1, \dots, g_w^r\colon A^m\to A^n$
defined by
\[
g^i_w(x) = x'\circ w_i,
\]
where $x\in A^m$. Formulas \eqref{E:isvs} imply that this is a strong sequence.

The following result is a version of the Halpern--L\"{a}uchli theorem (for strong subtrees). Recall
the definition of the trees $T^{k,n}$ from Section~\ref{S:tree}. Fix $k$ and let $T^n = T^{k,n}$. Note that
we can take $T^n = A^{\leq n}$ for $A$ with $|A|=k$.

{\bf (HL1)} {\em Given
$d>0$, $t$ and $m$ there exists $n$ such that for each $d$-coloring
of ${\rm leaves}(T^n)\times \cdots \times {\rm leaves}(T^n)$ ($t$ factors) there exists a strong sequence of leaf preserving
embeddings $g_i\colon T^m\to T^n$, for $i=1, \dots, t$, such that the set
\[
g_1({\rm leaves}(T^m))\times \cdots \times g_t({\rm leaves}(T^m))
\]
is monochromatic.}

{\bf (HL2)}  {\em Given $d>0$, $t$ and $m$ there exists $n$ such that for each $d$-coloring of
\[
\{ (w_1, \dots, w_t)\colon w_1, \dots, w_r\in T^n,\, {\rm ht}(w_i) = {\rm ht}(w_j),\hbox{ for }1\leq i,j\leq t\}
\]
there exists a strong sequence of
embeddings $g_i\colon T^m\to T^n$ for $i=1, \dots, t$ such that the set
\[
\{ (g_1(v_1), \dots g_t(v_t))\colon v_1, \dots, v_t\in T^m,\, {\rm ht}(v_i) = {\rm ht}(v_j),\hbox{ for }
1\leq i,j\leq t\}
\]
is monochromatic.}

We show that the above statements are re-phrasings of the
Hales--Jewett theorem. The Hales--Jewett theorem can be
stated as below in points (a) and (b).
(It is stated this way in \cite[Section 7]{Sol}, and it is proved there using the abstract approach to Ramsey theory.)

(a) Let $B$ be a finite set not including any natural numbers. Given
$d>0$ and $m$ there is $n$ such that for each $d$-coloring of
functions from $[n]$ to $B$ there is a function $w_0\colon [n]\to
B\cup [m]$ with properties (i) and (ii) such that the set
\[
\{ v\circ w_0\colon v \colon B\cup [m]\to B,\, v\res B = {\rm
id}_B\}
\]
is monochromatic.

(b) Let $B$ be a finite set not including any natural numbers. Given
$d>0$ and $m$ there is $n$ such that for each $d$-coloring of
functions from $[q]$ to $B$ for all $q\leq n$ there is $n_0\leq n$
and a function $w_0\colon [n_0]\to B\cup [m]$ with properties (i) and (ii) such that the set
\[
\{ v\circ w_0\colon v\colon B\cup [p]\to B,\, p\leq m,\, v\res B =
{\rm id}_B\}
\]
is monochromatic.

By the discussion at the beginning of this appendix, it is clear that (HL1) and (HL2) follow from (a) and (b), respectively, by taking $B=A^t$.

\medskip

\noindent{\bf Acknowledgement.} I thank Kostya Slutsky for preparing
for me the two wonderful drawings, and Alekos Kechris, Justin Moore, Miodrag Soki{\'c}, Min Zhao, and 
the referee for their comments on earlier versions of the paper.


\begin{thebibliography}{10}
\bibitem{Cam} P.J. Cameron, {\em Some treelike objects}, Quart. J. Math. Oxford 38 (1987), 155--183.

\bibitem{Deu} W. Deuber, {\em A generalization of Ramsey's theorem
for regular trees}, J. Combin. Theory, Ser. B 18 (1975), 18--23.

\bibitem{Jas} J. Jasi{\'n}ski, {\em Properties of boron tree structures}, preprint, 2011.

\bibitem{Mil} K. Milliken, {\em A Ramsey theorem for tress}, J.
Combin. Theory, Ser. A 26 (1979), 137--148.

\bibitem{Nes} J. Ne\v{s}et\v{r}il, {\em Ramsey theory}, in {\em
Handbook of Combinatorics}, eds. R. Graham, M. Gr\"{o}tschel, L.
Lov\'{a}sz, Elsevier Science, 1995, pp. 1331--1403.

\bibitem{Soc} M. Soki{\'c}, {\em Bounds on trees}, Discrete Math. 311 (2011), 398--407.

\bibitem{Sol} S. Solecki, {\em Abstract approach to finite Ramsey
theory and a self-dual Ramsey theorem}, preprint, 2011.
\end{thebibliography}
\end{document}